\documentclass[graybox]{svmult}

\usepackage{type1cm}        % activate if the above 3 fonts are
\usepackage{makeidx}         % allows index generation
\usepackage{graphicx}        % standard LaTeX graphics tool
\usepackage{tabu}
\usepackage[caption=false]{subfig}
\usepackage{tikz}
\usepackage{url}
\usepackage{multicol}        % used for the two-column index
\usepackage{multirow}
\usepackage[bottom]{footmisc}% places footnotes at page bottom

\usepackage{newtxtext}       %
\usepackage{newtxmath}       % selects Times Roman as basic font

\newcommand{\VDelta}{V_{\Delta}^0}
\newcommand{\VH}{V_{\mathcal H}^0}
\newcommand{\cRDelta}{\mathcal{R}_{0,\Delta}}
\newcommand{\cRH}{\mathcal{R}_{0,\mathcal H}}
\newcommand{\RDelta}{{R}_{0,\Delta}}
\newcommand{\RH}{{R}_{0,\mathcal H}}
\newcommand{\MDelta}{M_{\mathrm{AS}, \Delta}}
\newcommand{\MH}{M_{\mathrm{AS}, \mathcal H}}
\newcommand{\BDelta}{B_{0,\Delta}}
\newcommand{\BH}{B_{0, \mathcal H}}
\newcommand{\z}{\hphantom{0}}

\makeindex             % used for the subject index
\begin{document}

\title*{GenEO coarse spaces for heterogeneous indefinite elliptic problems}
\titlerunning{GenEO for heterogeneous indefinite elliptic problems}
\author{Niall Bootland, Victorita Dolean, Ivan G. Graham, Chupeng Ma and Robert Scheichl}
\authorrunning{N. Bootland et al.}
\institute{Niall Bootland \at University of Strathclyde, Dept. of Maths and Stats, \email{niall.bootland@strath.ac.uk}
\and Victorita Dolean \at University of Strathclyde, Dept. of Maths and Stats and University C\^ote d'Azur, CNRS, LJAD \email{work@victoritadolean.com}
\and Ivan G. Graham \at University of Bath, Dept. Math. Sci., \email{i.g.graham@bath.ac.uk}
\and Chupeng Ma \at Heidelberg University, Inst. Applied Math., \email{chupeng.ma@uni-heidelberg.de}
\and Robert Scheichl \at Heidelberg University, Inst. Applied Math., \email{r.scheichl@uni-heidelberg.de}}
%
% Use the package "url.sty" to avoid
% problems with special characters
% used in your e-mail or web address
%
\maketitle

%\abstract*{
\abstract{
Motivated by recent work on coarse spaces for Helmholtz problems, we provide in this paper a comparative study on the use of spectral coarse spaces of GenEO type for heterogeneous indefinite elliptic problems within an additive overlapping Schwarz method. In particular, we focus here on two different but related formulations of local generalised eigenvalue problems and compare their performance numerically. Even though their behaviour seems to be very similar for several well-known heterogeneous test cases that are mildly indefinite, only one of the coarse spaces has so far been analysed theoretically, while the other one leads to a significantly more robust domain decomposition method when the indefiniteness is increased. We present a summary of upcoming results developing such a theory and describe how the numerical experiments illustrate it.
}

\section{Introduction and motivations}
\label{sec:intro}

For domain decomposition preconditioners, the use of a coarse correction as a second level is usually required to provide scalability (in the weak sense), such that the iteration count is independent of the number of subdomains, for subdomains of fixed dimension. In addition, it is desirable to guarantee robustness with respect to strong variations in the physical parameters. Achieving scalability and robustness usually relies on sophisticated tools such as spectral coarse spaces \cite{Galvis:2010,Dolean:2015:DDM}. In particular, we can highlight the GenEO coarse space \cite{Spillane:2014:ARC}, which has been successfully analysed and applied to highly heterogeneous positive definite elliptic problems. This coarse space relies on the solution of local eigenvalue problems on subdomains and the theory in the SPD case is based on the fact that local eigenfunctions form an orthonormal basis with respect to the energy scalar product induced by the bilinear form.

Our motivation here is to gain a better insight into the good performance of spectral coarse spaces even for highly indefinite high-frequency Helmholtz problems with absorbing boundary conditions, as observed in \cite{Conen:2014:ACS} (for the Dirichlet-to-Neumann coarse space) and more recently in \cite{Bootland:2020:CSC} for coarse spaces of GenEO type. While a rigorous analysis for Helmholtz problems still lies beyond reach (see also \cite{Gander:2011} for the challenges), we present here numerical results, showing the benefits of GenEO-type coarse spaces for the heterogeneous symmetric indefinite elliptic problem
\begin{align}
\label{HelmholtzEquation}
- \nabla \cdot (A(\vec{x}) \nabla u) - \kappa u = f \ \ \text{in}\ \ \Omega\, , \quad \quad \text{subject to} \quad \quad u = 0 \ \ \text{on}\ \ \partial\Omega\, ,
\end{align}
in a bounded domain $\Omega$ with homogeneous Dirichlet boundary conditions on $\partial\Omega$, thus extending the results of \cite{Spillane:2014:ARC} to this case. The coefficient function $A$ in \eqref{HelmholtzEquation} is a symmetric positive-definite matrix-valued function on $\Omega \rightarrow {\mathbb R}^{d\times d}$ (where $d$ is the space dimension) with highly varying but bounded values ($a_{\rm min} |\xi|^{2}\leq {A}(\vec{x})\xi\cdot\xi \leq {a_{\rm max} |\xi|^{2}, \vec{x} \in\Omega, \xi \in\mathbb{R}^{d}}$) and $\kappa$ is an $L^\infty(\Omega)$ function which can have positive or negative values. We assume throughout that problem \eqref{HelmholtzEquation} is well-posed and that there is a unique weak solution $u \in H^{1}_{0}(\Omega)$, for all $f \in L^2(\Omega)$.

We propose two types of spectral coarse spaces, one built from local spectra of the whole indefinite operator on the left-hand side of \eqref{HelmholtzEquation}, and the other built using only the second-order operator in \eqref{HelmholtzEquation}. For the latter, the analysis in \cite{Bootland:2021:OSM} will apply, while the better performance of the former for large $\Vert \kappa \Vert_{\infty}$ provides some insight into the good performance of the $\mathcal{H}$-GenEO method introduced in \cite{Bootland:2020:CSC} for high-frequency Helmholtz problems, even though it is not amenable to the theory in~\cite{Bootland:2021:OSM}.

The problem \eqref{HelmholtzEquation} involves a Helmholtz-type operator (although this term would normally be associated with the case when $\kappa$ has a positive sign and \eqref{HelmholtzEquation} would normally be equipped with an absorbing boundary condition rather than the Dirichlet condition here). In the special case $A = I, \kappa = k^2$ with $k $ constant, the assumption of well-posedness of the problem is equivalent to the requirement that $k^2$ does not coincide with any of the Dirichlet eigenvalues of the operator $- \Delta $ in the domain $\Omega$. In this case, for large $k^2$, the solution of \eqref{HelmholtzEquation} will be rich in modes corresponding to eigenvalues near $k^2$ and thus will have oscillatory behaviour, increasing as $k$ increases. The Helmholtz problem with $A = I$ and $\kappa = \omega^2 n$ (with $\omega$ real and $n$ a function), together with an absorbing far-field boundary condition appears regularly in geophysical applications; here $n$ is the refractive index or `squared slowness' of waves and $\omega$ is the angular frequency.

To solve discretisations of \eqref{HelmholtzEquation}, we consider an additive Schwarz (AS) method with a GenEO-like coarse space and study the performance of this solver methodology for some heterogeneous test cases. GenEO coarse spaces have been shown theoretically and practically to be very effective for heterogeneous positive definite problems. Here, our main focus is to investigate how this approach performs in the indefinite case~\eqref{HelmholtzEquation}. We now review the underlying numerical methods that are used.

\section{Discretisation and domain decomposition solver}
\label{sec:discretisation}

We suppose that the domain $\Omega$ is a bounded Lipschitz polygon/polyhedron in 2D/3D. To discretise the problem we use the Lagrange finite element method of degree $p$ on a conforming simplicial mesh $T^{h}$ of $\Omega$. Denote the finite element space by $V^{h} \subset H^{1}_0(\Omega)$. The finite element solution $u_{h} \in V^{h}$ satisfies the weak formulation $b(u_{h},v_{h}) = F(v_{h})$, for all $v_{h} \in V^{h}$, where
\begin{align}
\label{WeakForms}
b(u,v) & = \int_{\Omega} \left( A(\vec{x})\nabla u \cdot \nabla {v} - \kappa u {v}\right) \, \mathrm{d}\vec{x} & \text{and} & & F(v) & = \int_{\Omega} {f v} \, \mathrm{d}\vec{x}.
\end{align}
Using the standard nodal basis for $V^{h}$ we can represent the solution $u_{h}$ through its basis coefficients $\vec{u}$ and reduce the problem to solving the symmetric linear system
\begin{align}\label{global}
B\vec{u} = \vec{f}
\end{align}
where $B$ comes from the bilinear form $b(\cdot,\cdot)$ and $\vec{f}$ from the linear functional $F(\cdot)$. Note that $B$ is symmetric but generally indefinite. For sufficiently small fine-mesh diameter $h$, problem~\eqref{global} has a unique solution $\vec{u}$; see \cite{Schatz:1996:SNE}. To solve \eqref{global}, we utilise a two-level domain decomposition preconditioner within a Krylov method.

Consider an overlapping partition $\lbrace \Omega_{j} \rbrace_{1 \le j \le N}$ of $\Omega$, where each $\Omega_{j}$ is assumed to have diameter $H_{j}$ and $H$ denotes the maximal diameter of the subdomains. For each $j$ we define $\widetilde{V}_{j} = \{v|_{\Omega_{j}} : v \in V^{h}\}$, $V_{j} = \{v \in \widetilde{V}_{j} : {\rm supp}(v)\subset \Omega_{j}\}$, and for $u$, $v \in \widetilde{V}_{j}$
\begin{align*}
b_j(u,v) &:= \int_{\Omega_{j}} \left(A(\vec{x})\nabla u\cdot{\nabla {v}} - \kappa u {v}\right) \mathrm{d}\vec{x} & \text{and} & & a_j(u,v) &:= \int_{\Omega_{j}} A(\vec{x}) \nabla u\cdot{\nabla {v}} \mathrm{d}\vec{x}.
\end{align*}

Let $ {\cal R}_{j}^{T} \colon V_{j} \rightarrow V^{h}$, $1\le j \le N$, denote the zero-extension operator, let $R_j^T$ denote its matrix representation with respect to the nodal basis and set $R_j = (R_j^T)^T$. The classical one-level additive Schwarz preconditioner is
\begin{align}\label{eq:2-14}
M_{\text{AS}}^{-1} &= \sum_{j=1}^{N} R_{j}^{T} B_{j}^{-1} R_{j}, & \text{where} & & B_{j} &= R_{j} B R_{j}^{T}.
\end{align}
It is well-known that one-level additive Schwarz methods are not scalable with respect to the number of subdomains in general, since information is exchanged only between neighbouring subdomains. Thus, we introduce the two-level additive Schwarz method with GenEO coarse space first proposed in \cite{Spillane:2014:ARC}. To this end, for $1\leq j\leq N$, let $\{\phi^{j}_{1},\ldots,\phi^{j}_{\tilde{n}_{j}}\}$ be a nodal basis of $\widetilde{V}_{j}$, where $\widetilde{n}_{j} = \mathrm{dim} ( \widetilde{V}_{j} )$.

\begin{definition}[Partition of unity] \label{def:2-0}
Let $\text{dof}(\Omega_j)$ denote the internal degrees of freedom (nodes) on subdomain $\Omega_j$. For any degree of freedom $i$, let $\mu_{i}$ denote the number of subdomains $\Omega_j$ for which $i$ is an internal degree of freedom, i.e., $\mu_{i} := \#\{j : 1 \leq j \leq N,\ i \in \text{dof}(\Omega_j)\}$. Then, for $1\leq j\leq N$, the \emph{local partition of unity operator} $\Xi_{j} \colon \widetilde{V}_{j}\rightarrow V_{j}$ is defined by
\begin{align}\label{eq:2-16}
\Xi_{j}(v) &:= \sum_{i \in \text{dof}(\Omega_j)}\frac{1}{\mu_i} v_{i} \phi^{j}_{i}, & \text{for all} & & v &= \sum_{i=1}^{\widetilde{n}_{j}}v_{i}\phi^{j}_{i} \in \widetilde{V}_{j}.
\end{align}
The operators $\Xi_{j}$ form a partition of unity, i.e., $\sum_{j=1}^{N} R_{j}^{T} \Xi_{j}(v|_{\Omega_j}) = v$, $\forall v \in V^{h}$ \cite{Spillane:2014:ARC}.
\end{definition}

\noindent
For each $j$, we define the following generalised eigenvalue problems:
\begin{align}
\label{eq:deltag}
\text{find} \ p \in \widetilde{V}_j\backslash \{0\} , \, \lambda \in \mathbb{R} : & & a_{j}(p,v) &= \lambda \,a_{j}(\Xi_{j}(p), \Xi_{j}(v)), & & \text{for all } v \in \widetilde{V}_{j}, \\
\label{eq:hg}
\text{find} \ q \in \widetilde{V}_j\backslash \{0\} , \, \lambda \in \mathbb{R} : & & b_{j}(q,v) &= \lambda \,a_{j}(\Xi_{j}(q), \Xi_{j}(v)), & & \text{for all } v \in \widetilde{V}_{j},
\end{align}
where $\Xi_{j}$ is the local partition of unity operator from Definition~\ref{def:2-0}.

\begin{definition}[$\Delta$-GenEO and ${\cal H}$-GenEO coarse spaces]
\label{def:2-2}
For each $j$, $1 \le j \le N$, let $(p_{l}^{j})_{l=1}^{m_{j}}$ and $(q_{l}^{j})_{l=1}^{m_{j}}$ be the eigenfunctions of the eigenproblems \eqref{eq:deltag} and \eqref{eq:hg} corresponding to the $m_{j}$ smallest eigenvalues, respectively. Then we define the \emph{$\Delta$-GenEO} and \emph{${\cal H}$-GenEO} coarse spaces, respectively, by
\begin{align}
\label{eq:deltspace}
\VDelta &:= {\rm span}\{{\cal R}_{j}^{T}\Xi_{j}(p_{l}^{j}): \;l=1,\ldots,m_{j};\;j=1,\ldots,N\} \ \ \text{and}\\
\label{eq:hspace}
\VH &:= {\rm span}\{{\cal R}_{j}^{T}\Xi_{j}(q_{l}^{j}): \;l=1,\ldots,m_{j};\;j=1,\ldots,N\}.
\end{align}
\end{definition}

Note that here and subsequently, the subscript $\Delta$ refers to the GenEO coarse space \eqref{eq:deltspace} based on \eqref{eq:deltag}, the eigenproblem with respect to the `Laplace-like' operator induced by the bilinear form $a_j$, while the subscript ${\cal H}$ refers to the ${\cal H}$-GenEO coarse space \eqref{eq:hspace} based on \eqref{eq:hg}, with the `Helmholtz-like' operator appearing in $b_j$.

Since $\VDelta, \, \VH \subset V^{h}$, we can introduce the natural embeddings $\cRDelta^T \colon \VDelta \rightarrow V^h$ and $\cRH^T \colon \VH \rightarrow V^h$, with matrix representations $\RDelta^T$ and $\RH^T$, respectively, and set $\RDelta = (\RDelta^T)^T$ and $\RH = (\RH^T)^T$ to obtain the following two-level extensions of the one-level additive Schwarz method \eqref{eq:2-14}:
\begin{align}\label{twolevel}
\MDelta^{-1} &= M_{AS}^{-1} + \RDelta^{T}\BDelta^{-1}\RDelta & \text{and} & & \MH^{-1} &= M_{AS}^{-1} + \RH^{T}\BH^{-1}\RH,
\end{align}
where $\BDelta := \RDelta {B} \RDelta^{T}$ and $\BH := \RH{B} \RH^{T}$.

\section{Theoretical results}
\label{sec:theory}

The theoretical properties of the preconditioner $\MDelta^{-1}$ are studied in the forthcoming paper \cite{Bootland:2021:OSM}. There, the PDE studied is a generalisation of \eqref{HelmholtzEquation}, which also allows the inclusion of a non-self-adjoint first order convection term. The important parameters in the preconditioner are the coarse mesh diameter $H$ and the `eigenvalue tolerance'
\begin{align*}
\Theta := \max_{1 \le j \le N } \left( \lambda_{m_j+1}^j \right)^{-1},
\end{align*}
where $\{\lambda_{m}^j: m = 1, 2, \ldots\}$ are the eigenvalues of the generalised eigenproblem~\eqref{eq:deltag}, given in non-decreasing order. We now highlight a special case of the results in \cite{Bootland:2021:OSM}.
\begin{theorem}
\label{thm:main}
Let the fine-mesh diameter $h$ be sufficiently small. Then there exist thresholds $H_0>0$ and $\Theta_0 > 0 $ such that, for all $H \leq H_0$ and $\Theta \leq \Theta_0$: the matrices $B_j$ and $\BDelta$ appearing in \eqref{eq:2-14} and \eqref{twolevel} are non-singular. Moreover, if problem \eqref{global} is solved by GMRES with left preconditioner $\MDelta^{-1}$ and residual minimisation in the energy norm $\Vert u \Vert_a := \big(\int_\Omega \nabla u \cdot A \nabla u \big)^{1/2} $, then there exists a constant $c \in (0,1)$, which depends on $H_0$ and $\Theta_0$ but is independent of all other parameters, such that we have the robust GMRES convergence estimate
\begin{equation}\label{eq:2-31}
\Vert r_{\ell} \Vert^{2}_{a} \leq \left(1-c^2\right)^{\ell}\Vert r_{0} \Vert^{2}_{a}\,,
\end{equation}
for $\ell = 0, 1, \ldots \,$, where $r_\ell$ denotes the residual after $\ell$ iterations of GMRES.
\end{theorem}
In fact, the paper \cite{Bootland:2021:OSM} will investigate in detail how the thresholds $H_0$ and $\Theta_0$ depend on the heterogeneity and indefiniteness of \eqref{HelmholtzEquation}. For example, if the problem is scaled so that $a_{\min} = 1$, then as $\Vert \kappa \Vert_\infty$ grows, $H_0$ and $\Theta_0$ have to decrease to maintain the convergence rate of GMRES:
\begin{align}
H_0 &\lesssim \Vert \kappa \Vert_\infty^{-1} & \text{and} & & \Theta_0 &\lesssim C_{\mathrm{stab}}^{-2} \Vert \kappa \Vert_\infty^{-4} \,, \label{rates}
\end{align}
where $C_{\mathrm{stab}} = C_{\mathrm{stab}}(A,\kappa)$ denotes the stability constant for problem \eqref{HelmholtzEquation}, i.e., the solution $u$ satisfies $\Vert u \Vert_{H^1(\Omega)} \leq C_{\mathrm{stab}} \Vert f \Vert_{L^2(\Omega)}$ for all $f \in L^2(\Omega)$ and the hidden constants are independent of $h$, $H$, $a_{\max}$ and $\kappa$. Thus, as $\Vert \kappa \Vert_\infty$ gets smaller, the indefiniteness diminishes and the requirements on $H_0$ and $\Theta_0$ are relaxed.

\section{Numerical results}
\label{sec:numerics}

We give results for a more efficient variant of the preconditioner described in \S\ref{sec:discretisation}. Instead of \eqref{twolevel}, we here use the \emph{restricted additive Schwarz} (RAS) method, with the GenEO coarse space incorporated using a deflation approach, yielding:
\begin{align}
\label{2LevelAdaptiveDeflationPreconditioner}
M^{-1} &= M_{\text{RAS}}^{-1}(I - B Q_0) + Q_0, & \quad \text{where} \quad M_{\text{RAS}}^{-1} &= \sum_{j=1}^{N} R_{j}^{T} D_j B_{j}^{-1} R_{j}.
\end{align}
Here, $D_j$ is the matrix form of the partition of unity operator $\Xi_{j}$. Moreover, we have $Q_0 = R_{0}^{T}B_{0}^{-1}R_{0}$ with $B_{0} = R_{0}BR_{0}^{T}$ and either $R_0 = \RDelta$ or $R_0 = \RH$, depending on whether we use $\Delta$-GenEO or $\mathcal{H}$-GenEO. We include all eigenfunctions $p_l^j$ or $q_l^j$ in $\VDelta$ or $\VH$ corresponding to eigenvalues $\lambda_l^j < \lambda_\text{max}$, for $\Delta$-GenEO or $\mathcal{H}$-GenEO, respectively. In $\mathcal{H}$-GenEO this includes all eigenfunctions corresponding to negative eigenvalues. Unless otherwise stated, the eigenvalue threshold is $\lambda_\text{max} = \frac{1}{2}$.

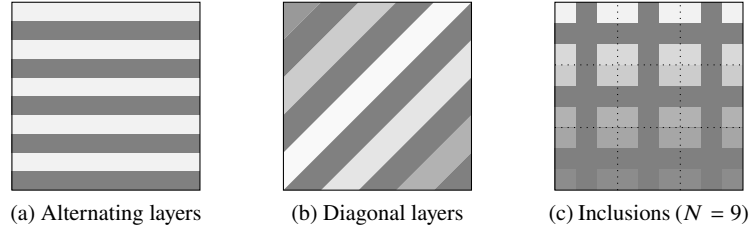
\begin{figure}[t]
	\centering
	\hspace*{\fill}
	\subfloat[Alternating layers]{
	\begin{tikzpicture}[scale=1.25]
	\fill[gray,opacity=0.1] (-1,1) -- (1,1) -- (1,0.8) -- (-1,0.8) -- cycle;
	\fill[gray,opacity=1.0] (-1,0.6) -- (1,0.6) -- (1,0.8) -- (-1,0.8) -- cycle;
	\fill[gray,opacity=0.1] (-1,0.6) -- (1,0.6) -- (1,0.4) -- (-1,0.4) -- cycle;
	\fill[gray,opacity=1.0] (-1,0.2) -- (1,0.2) -- (1,0.4) -- (-1,0.4) -- cycle;
	\fill[gray,opacity=0.1] (-1,0.2) -- (1,0.2) -- (1,0.0) -- (-1,0.0) -- cycle;
	\fill[gray,opacity=1.0] (-1,-0.2) -- (1,-0.2) -- (1,0.0) -- (-1,0.0) -- cycle;
	\fill[gray,opacity=0.1] (-1,-0.2) -- (1,-0.2) -- (1,-0.4) -- (-1,-0.4) -- cycle;
	\fill[gray,opacity=1.0] (-1,-0.6) -- (1,-0.6) -- (1,-0.4) -- (-1,-0.4) -- cycle;
	\fill[gray,opacity=0.1] (-1,-0.6) -- (1,-0.6) -- (1,-0.8) -- (-1,-0.8) -- cycle;
	\fill[gray,opacity=1.0] (-1,-1) -- (1,-1) -- (1,-0.8) -- (-1,-0.8) -- cycle;
	\draw (-1,1) -- (-1,-1) -- (1,-1) -- (1,1) -- cycle;
	\label{Fig:LayersB}
	\end{tikzpicture}
	} \hspace*{\fill}
	\subfloat[Diagonal layers]{
	\begin{tikzpicture}[scale=1.25]
	\fill[gray,opacity=1.0] (1,-1) -- (0.6,-1) -- (1,-0.6) -- cycle;
	\fill[gray,opacity=0.6] (0.6,-1) -- (0.2,-1) -- (1,-0.2) -- (1,-0.6) -- cycle;
	\fill[gray,opacity=1.0] (0.2,-1) -- (-0.2,-1) -- (1,0.2) -- (1,-0.2) -- cycle;
	\fill[gray,opacity=0.2] (-0.2,-1) -- (-0.6,-1) -- (1,0.6) -- (1,0.2) -- cycle;
	\fill[gray,opacity=1.0] (-0.6,-1) -- (-1,-1) -- (1,1) -- (1,0.6) -- cycle;
	\fill[gray,opacity=0.05] (-1,-1) -- (-1,-0.6) -- (0.6,1) -- (1,1) -- cycle;
	\fill[gray,opacity=1.0] (-1,-0.6) -- (-1,-0.2) -- (0.2,1) -- (0.6,1) -- cycle;
	\fill[gray,opacity=0.4] (-1,-0.2) -- (-1,0.2) -- (-0.2,1) -- (0.2,1) -- cycle;
	\fill[gray,opacity=1.0] (-1,0.2) -- (-1,0.6) -- (-0.6,1) -- (-0.2,1) -- cycle;
	\fill[gray,opacity=0.8] (-1,0.6) -- (-1,1) -- (-0.6,1) -- cycle;
	\draw (-1,1) -- (-1,-1) -- (1,-1) -- (1,1) -- cycle;
	\label{Fig:LayersC}
	\end{tikzpicture}
	}
	\hspace*{\fill}
	\subfloat[Inclusions ($N=9$)]{
	\begin{tikzpicture}[scale=1.25]
	\fill[gray,opacity=0.1] (-1,1) -- (1,1) -- (1,7/9) -- (-1,7/9) -- cycle;
	\fill[gray,opacity=0.3] (-1,5/9) -- (1,5/9) -- (1,3/9) -- (-1,3/9) -- cycle;
	\fill[gray,opacity=0.4] (-1,3/9) -- (1,3/9) -- (1,1/9) -- (-1,1/9) -- cycle;
	\fill[gray,opacity=0.6] (-1,-1/9) -- (1,-1/9) -- (1,-3/9) -- (-1,-3/9) -- cycle;
	\fill[gray,opacity=0.7] (-1,-3/9) -- (1,-3/9) -- (1,-5/9) -- (-1,-5/9) -- cycle;
	\fill[gray,opacity=0.9] (-1,-7/9) -- (1,-7/9) -- (1,-1) -- (-1,-1) -- cycle;
	\fill[gray,opacity=1.0] (-1,5/9) -- (1,5/9) -- (1,7/9) -- (-1,7/9) -- cycle;
	\fill[gray,opacity=1.0] (-1,1/9) -- (1,1/9) -- (1,-1/9) -- (-1,-1/9) -- cycle;
	\fill[gray,opacity=1.0] (-1,-5/9) -- (1,-5/9) -- (1,-7/9) -- (-1,-7/9) -- cycle;
	\fill[gray,opacity=1.0] (5/9,-1) -- (5/9,1) -- (7/9,1) -- (7/9,-1) -- cycle;
	\fill[gray,opacity=1.0] (1/9,-1) -- (1/9,1) -- (-1/9,1) -- (-1/9,-1) -- cycle;
	\fill[gray,opacity=1.0] (-5/9,-1) -- (-5/9,1) -- (-7/9,1) -- (-7/9,-1) -- cycle;
	\draw (-1,1) -- (-1,-1) -- (1,-1) -- (1,1) -- cycle;
	\draw[dotted] (-1,1/3) -- (1,1/3); \draw[dotted] (-1,-1/3) -- (1,-1/3);
	\draw[dotted] (1/3,-1) -- (1/3,1); \draw[dotted] (-1/3,-1) -- (-1/3,1);
	\label{Fig:Inclusions}
	\end{tikzpicture}
	} \hspace*{\fill}
	\caption{Piecewise constant profiles $a(\vec{x})$, where $A(\vec{x}) = a(\vec{x})I$. For the darkest shade $a(\mathbf{x}) = 1$ while for the lightest shade $a(\mathbf{x}) = a_\text{max}$. Profiles (a) and (b) are fixed while in (c) the interfaces of $a$ are the same per subdomain, although the value of $a$ depends on height (case $N = 3^2$ is shown).}
	\label{Fig:Layers}
\end{figure}

As a model problem, we consider \eqref{HelmholtzEquation} on the unit square $\Omega = (0,1)^2$, take $\kappa$ constant, and define $A$ to model various layered media, as depicted in Fig.~\ref{Fig:Layers}. The right-hand side $f$ is taken to be a point source at the centre $(\frac{1}{2},\frac{1}{2})$. To discretise, we use a uniform square grid with $n_{\text{glob}}$ points in each direction and triangulate along one set of diagonals to form P1 elements. We further use a uniform decomposition into $N$ square subdomains and throughout use minimal overlap (non-overlapping subdomains are extended by adding only the fine-mesh elements which touch them).

Our computations are performed using FreeFem (\url{http://freefem.org/}), in particular using the \texttt{ffddm} framework. We use preconditioned GMRES with residual minimisation in the Euclidean norm and a relative residual tolerance of $10^{-6}$. We have assumed $a_\text{min} = 1$; otherwise a rescaling will ensure this. The indefiniteness is controlled by $\kappa$, taken here to be a positive constant. Although estimate \eqref{eq:2-31} describes GMRES implemented in the energy inner product, we use here the standard Euclidean implementation and prove in \cite[\S 4]{Bootland:2021:OSM} that (for quasi-uniform meshes) the latter algorithm requires at most $\mathcal{O}(\log(a_{\text{max}}/h))$ more iterations than the former to achieve the same residual reduction. Experiments for Helmholtz problems in \cite{GrSpVa:17} showed that the two approaches performed almost identically.

In Table \ref{Table:VaryAmax_AlternatingLayers_nglob400} we provide GMRES iteration counts for $\Delta$-GenEO and $\mathcal{H}$-GenEO as $N$ varies in two cases: in case (i), on the left, we use profile (a) and increase $a_\text{max}$ while in case (ii), on the right, we use profile (c) and increase $n_\text{glob} = h^{-1}$. In (i) we see clear robustness to increasing the contrast parameter $a_\text{max}$. In (ii) we observe robustness to decreasing $h$, with markedly better performance for $\mathcal{H}$-GenEO. In (ii), the coefficient $a(\vec{x})$ (and hence the problem itself) becomes more complicated as $N$ increases since the geometry of the coefficient remains identical in each subdomain.

In Table \ref{Table:VaryC_DiagonalLayers_nglob600} we illustrate the effect of increasing $\kappa$, giving iteration counts and (in brackets) coarse space sizes. Here we see the substantial advantage of $\mathcal{H}$-GenEO over $\Delta$-GenEO: much better iteration counts are obtained, yet the coarse space size increases only modestly. As $N$ increases, although the dimension of the coarse space grows, the number of eigenfunctions per subdomain decreases. For very large $\kappa$, neither method is fully robust while, for small $\kappa$, both methods perform similarly. This leads to the interesting question of whether robustness to $\kappa$ can be gained by taking more eigenfunctions in the coarse space. Table \ref{Table:VaryLambdaMax_DiagonalLayers} gives results for a sequence of increasing values of $\kappa$ for the diagonal layers problem, in which we simultaneously increase $\lambda_{\text{max}}$, indicating (apparent almost) robustness with respect to $\kappa$.

\begin{table}[ht]
\centering
\caption{GMRES iteration counts with $\lambda_\text{max} = \frac{1}{2}$. Left-hand table: Alternating layers problem, varying $a_\text{max} > 1$ and $N$, with fixed $\kappa = 100$ and $n_\text{glob} = 400$. Right-hand table: Inclusions problem, varying $n_\text{glob}$ and $N$, with fixed $\kappa = 1000$ and $a_\text{max} = 50$.}
\label{Table:VaryAmax_AlternatingLayers_nglob400}
\tabulinesep=1mm
\begin{center}
\begin{tabular}{lr}
\begin{tabu}{|l|cccc|cccc|}\hline
	&\multicolumn{4}{c|}{$\Delta$-GenEO} & \multicolumn{4}{c|}{$\mathcal{H}$-GenEO} \\
	\cline{1-9}
	$a_\text{max} \ \backslash \ N$ & 16 & 36 & 64 & 100 & 16 & 36 & 64 & 100 \\
	\hline
	10   & 10 & 9 & 9 & 10  & 9 & 9 & 9 & 9 \\
	50   &  9 & 9 & 9 &  9  & 9 & 9 & 9 & 9 \\
	200  &  9 & 9 & 9 &  9  & 9 & 9 & 9 & 9 \\
	1000 &  9 & 9 & 9 &  9  & 9 & 9 & 9 & 9 \\
	\hline
\end{tabu}
\quad & \quad
\begin{tabu}{|l|cccc|cccc|}
	\hline & \multicolumn{4}{c|}{$\Delta$-GenEO} & \multicolumn{4}{c|}{$\mathcal{H}$-GenEO} \\
	\cline{1-9}
	$n_\text{glob} \ \backslash \ N$ & 16 & 36 & 64 & 100 & 16 & 36 & 64 & 100 \\
	\hline
	200 & 24 & 16 & 26 & 22  & 10 & 10 & 10 & 11 \\
	400 & 23 & 14 & 19 & 18  &  8 &  9 &  9 &  8 \\
	600 & 21 & 14 & 18 & 19  &  8 &  9 &  8 & 10 \\
	800 & 22 & 14 & 19 & 20  &  9 &  9 &  9 & 10 \\
	\hline
\end{tabu}
\end{tabular}
\end{center}
\end{table}

\begin{table}[ht]
\centering
\caption{GMRES iteration counts and (in brackets) coarse space dimension for the diagonal layers problem with $\lambda_\text{max} = \frac{1}{2}$, varying $\kappa$ and $N$, with fixed $a_\text{max} = 5$ and $n_\text{glob} = 600$.}
\label{Table:VaryC_DiagonalLayers_nglob600}
\tabulinesep=1mm
\begin{tabu}{|l|l|l|l|l|l|l|l|l|}
	\hline
	& \multicolumn{4}{c|}{$\Delta$-GenEO} & \multicolumn{4}{c|}{$\mathcal{H}$-GenEO} \\
	\cline{1-9}
	$\kappa \ \backslash \ N$ & 16 & 36 & 64 & 100 & 16 & 36 & 64 & 100 \\
	\hline
	10    & 9 \z\z(627) & 9 \z\z(1050) & 9 \z\z(1468) & 9 \z\z(1804)  & 9 \z\z(627) & 9 \z(1050) & 9 \z(1468) & 9 \z(1804) \\
	100   & 10  \z(627) & 9 \z\z(1050) & 9 \z\z(1468) & 9 \z\z(1804)  & 9 \z\z(627) & 9 \z(1052) & 9 \z(1473) & 9 \z(1814) \\
	1000  & 36  \z(627) & 43  \z(1050) & 35  \z(1468) & 28  \z(1804)  & 13  \z(674) & 11  (1083) & 9 \z(1520) & 10  (1877) \\
	10000 & 215   (627) & 339   (1050) & 437   (1468) & 506   (1804)  & 27   (1256) & 33  (1651) & 40  (2139) & 18  (2549) \\
	\hline
\end{tabu}
\end{table}

\begin{table}[ht]
\centering
\caption{GMRES iteration counts and (in brackets) coarse space dimension for $\mathcal{H}$-GenEO for the diagonal layers problem, varying $N$, with fixed $n_\text{glob} = 600$ and $a_\text{max} = 5$ and increasing eigenvalue threshold $\lambda_\text{max}$ as $\kappa$ increases, aiming to control iteration counts as $\kappa$ increases.}
\label{Table:VaryLambdaMax_DiagonalLayers}
\tabulinesep=1mm
\begin{tabu}{|ll|l|l|l|l|}
	\hline
	\multicolumn{2}{|c|}{} & \multicolumn{4}{c|}{$\mathcal{H}-$GenEO} \\
	\hline
	$\lambda_\text{max}$ & $\kappa \ \backslash \ N$ & 16 & 36 & 64 & 100 \\
	\hline
	0.1 & 10    & 23 \z(108) & 23 \z(199) & 25 \z(214) & 23 \z(324) \\
	0.1 & 100   & 23 \z(111) & 24 \z(201) & 28 \z(223) & 27 \z(324) \\
	0.2 & 1000  & 19 \z(265) & 27 \z(418) & 20 \z(574) & 20 \z(684) \\
	0.6 & 10000 & 24  (1430) & 25  (2129) & 28  (2680) & 15  (3252) \\
	\hline
\end{tabu}
\end{table}

\begin{figure}[t]
	\centering
	\tabulinesep=1mm
	\begin{tabu}{ccccc}
		& \multirow{2}{*}{$\Delta$-GenEO} & $\mathcal{H}$-GenEO & $\mathcal{H}$-GenEO & $\mathcal{H}$-GenEO \\
		& & $\kappa = 1000$ & $\kappa = 10000$ & $\kappa = 10000$ \\
		\rotatebox{90}{Homogeneous} &
		\includegraphics[width=.24\linewidth]{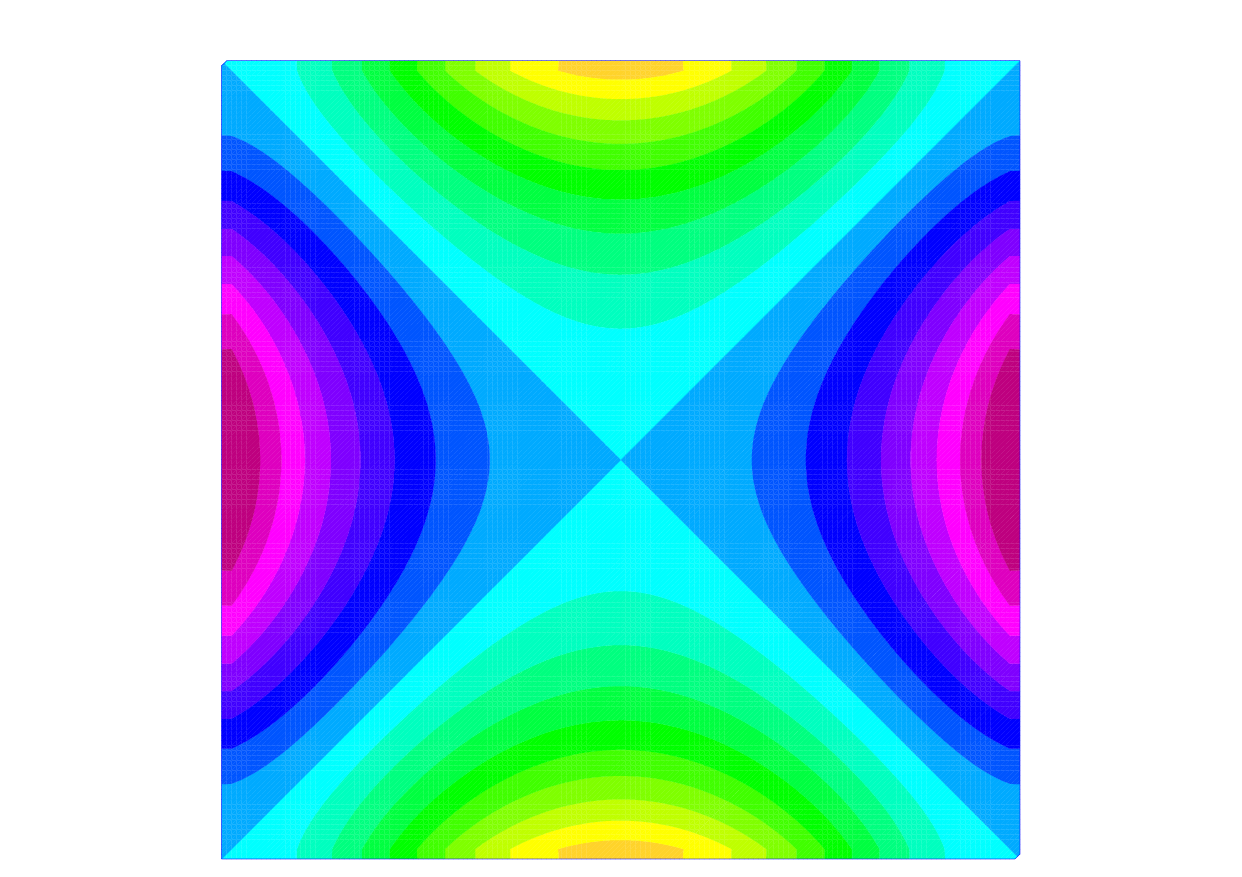} &
		\includegraphics[width=.24\linewidth]{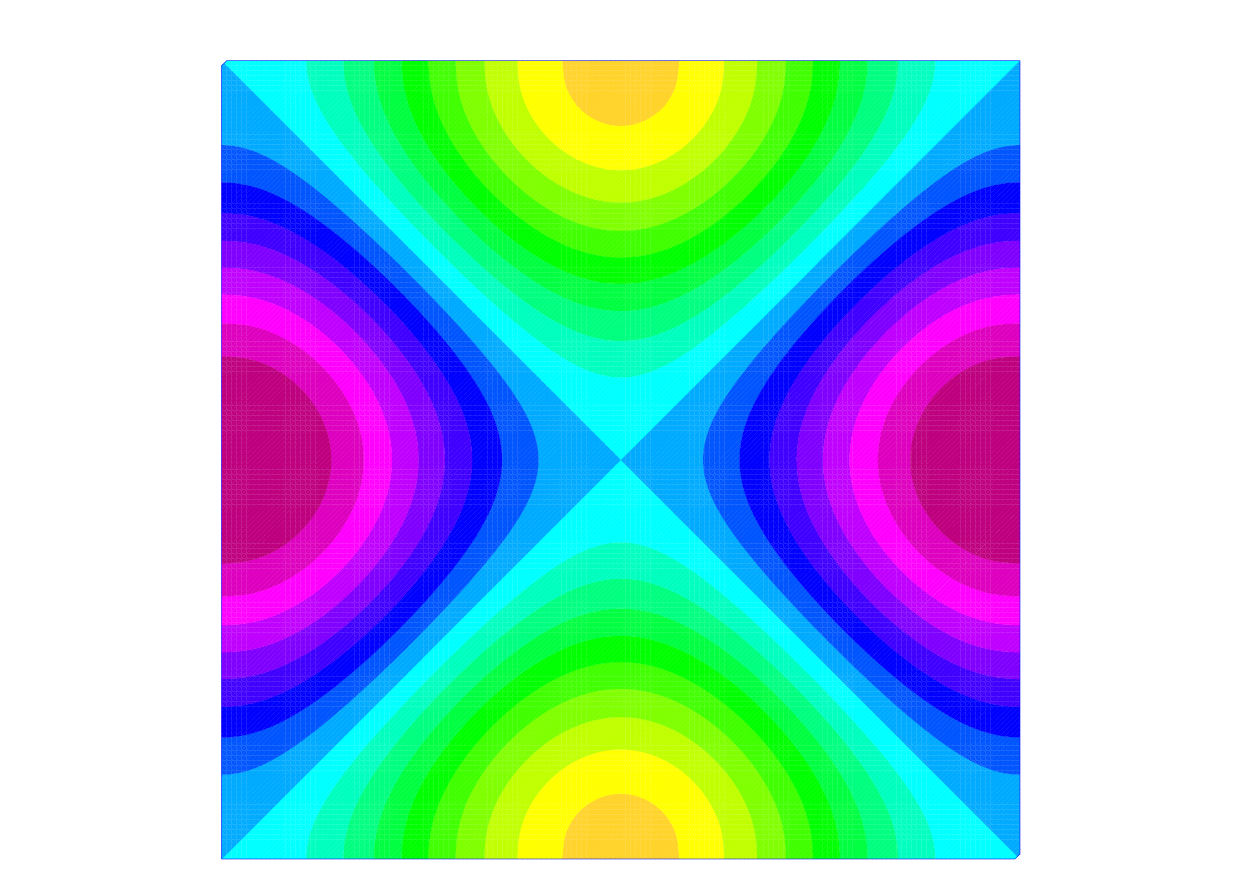} &
		\includegraphics[width=.24\linewidth]{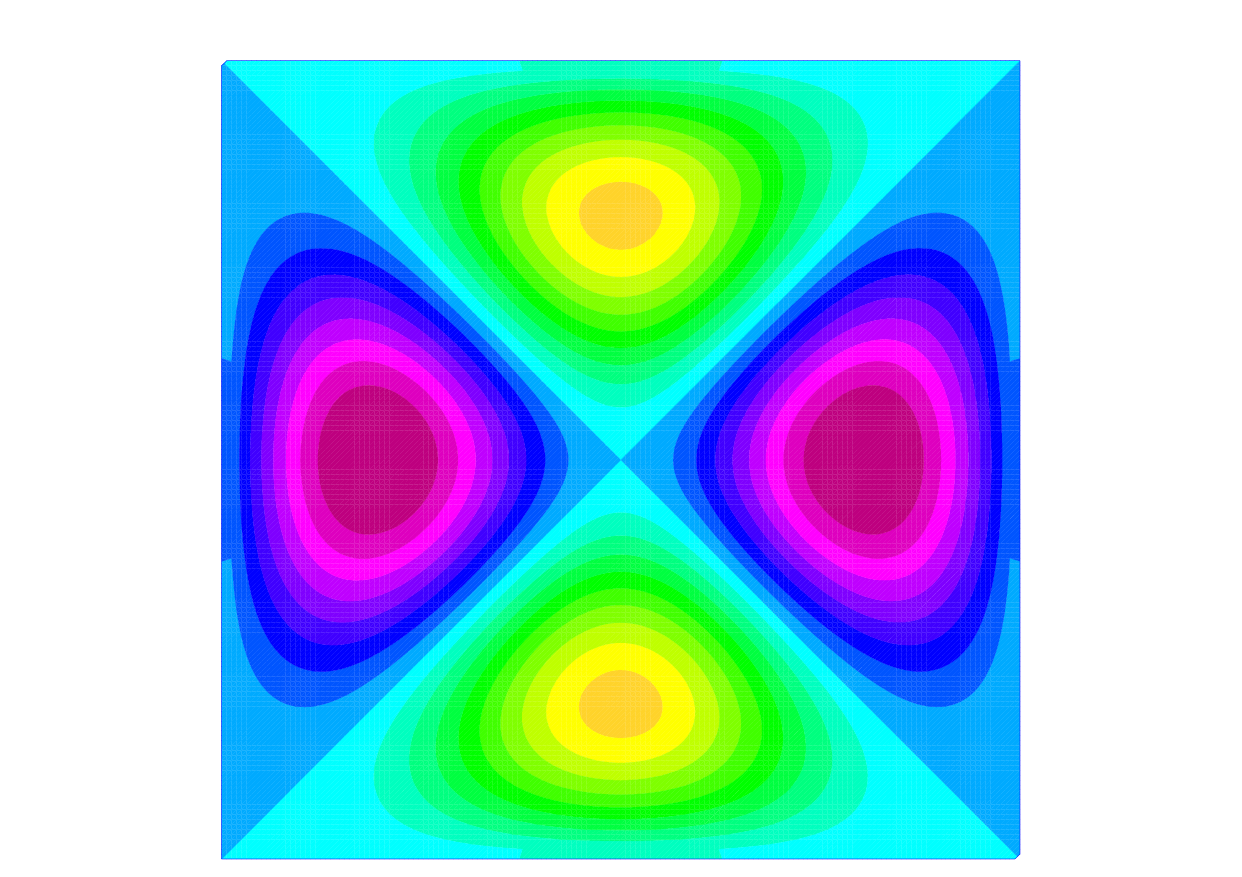} &
		\includegraphics[width=.24\linewidth]{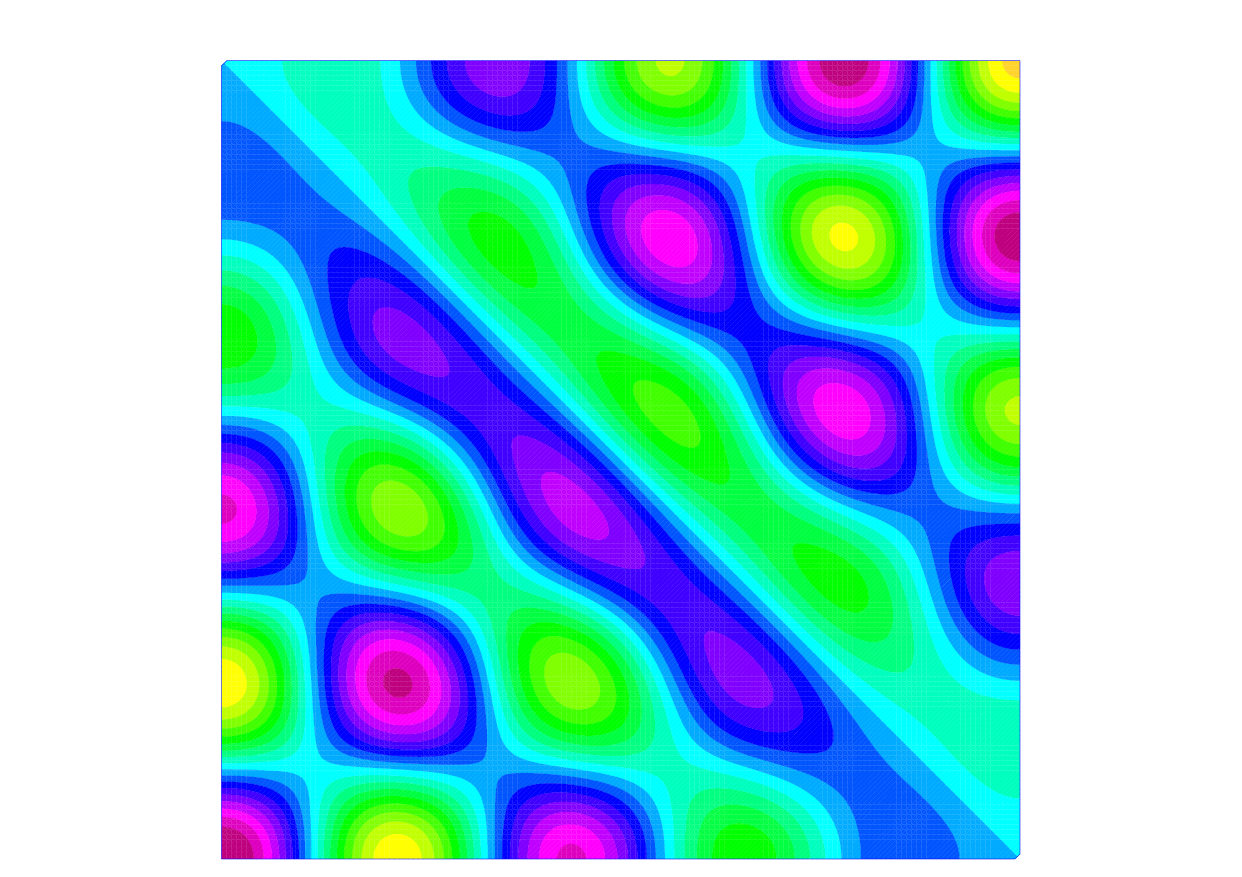} \\[-1ex]
		& $\lambda = 0.057$ & $\lambda = -0.003$ & $\lambda = -3.220$ & $\lambda = -0.004$ \\[1ex]
		\rotatebox{90}{Diagonal layers} &
		\includegraphics[width=.24\linewidth]{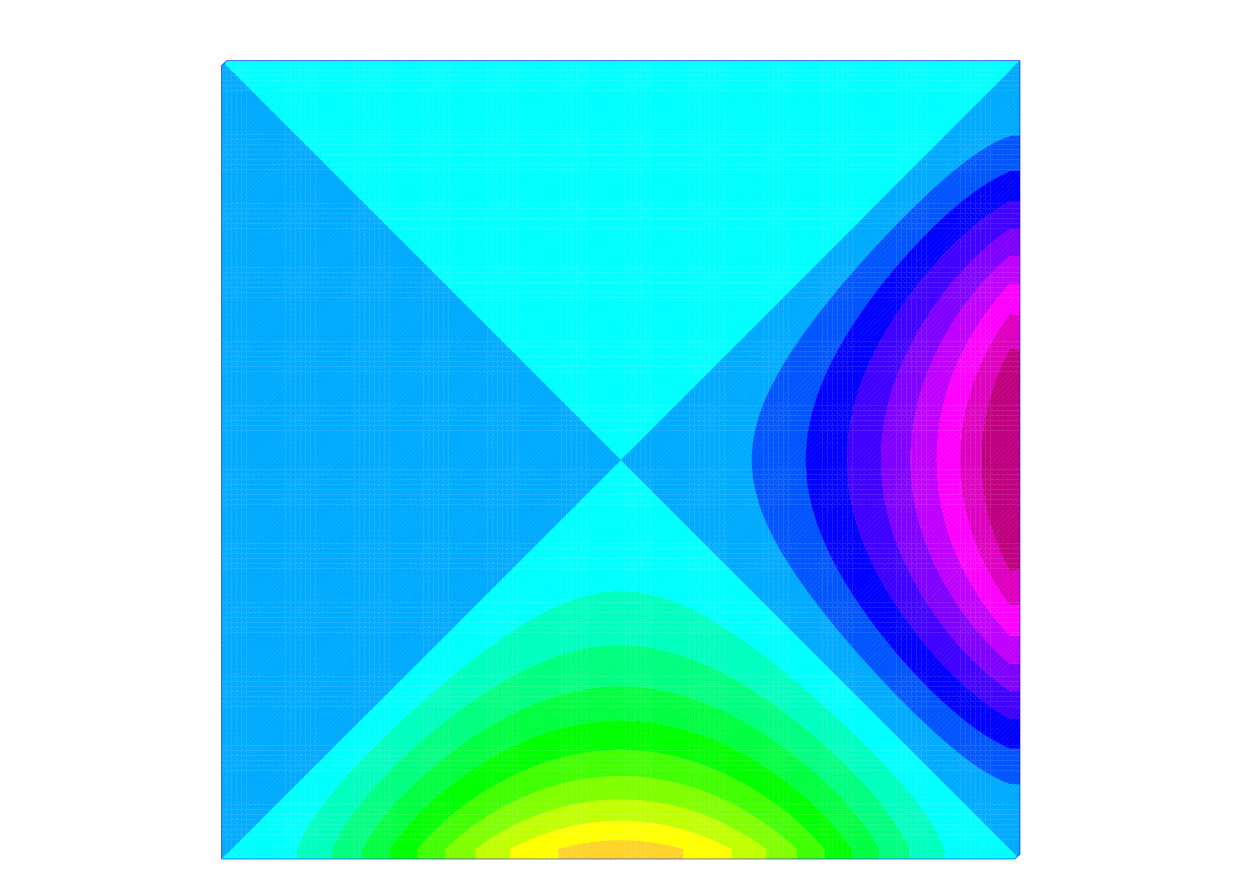} &
		\includegraphics[width=.24\linewidth]{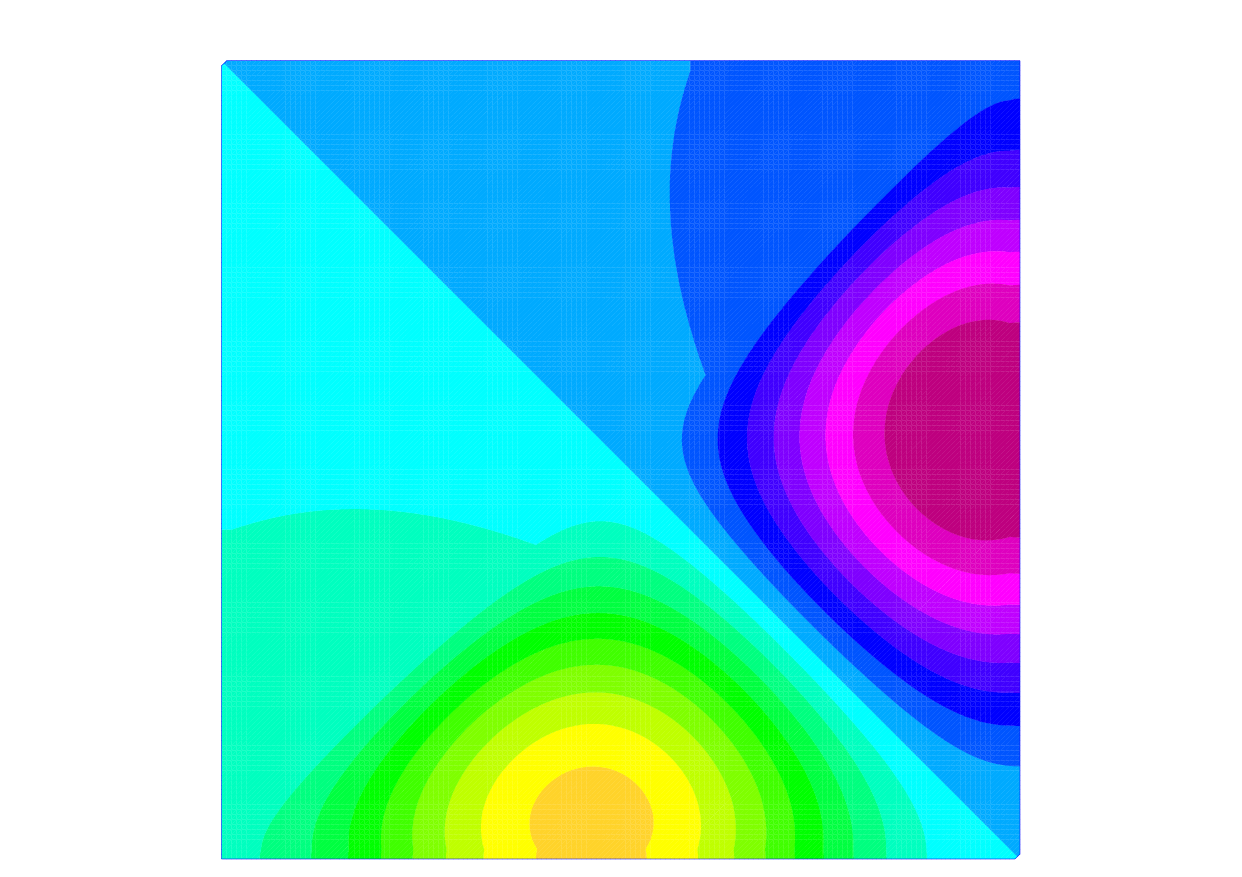} &
		\includegraphics[width=.24\linewidth]{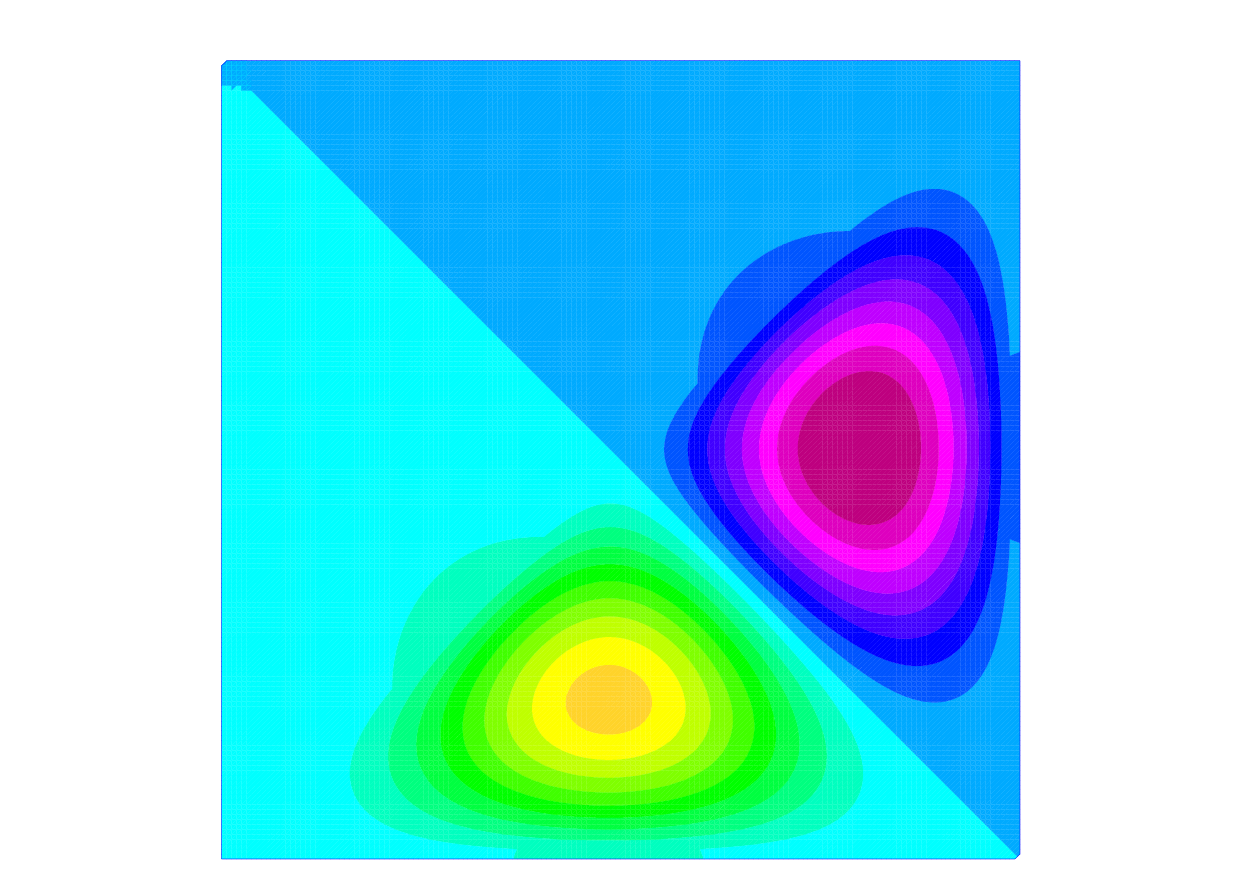} &
		\includegraphics[width=.24\linewidth]{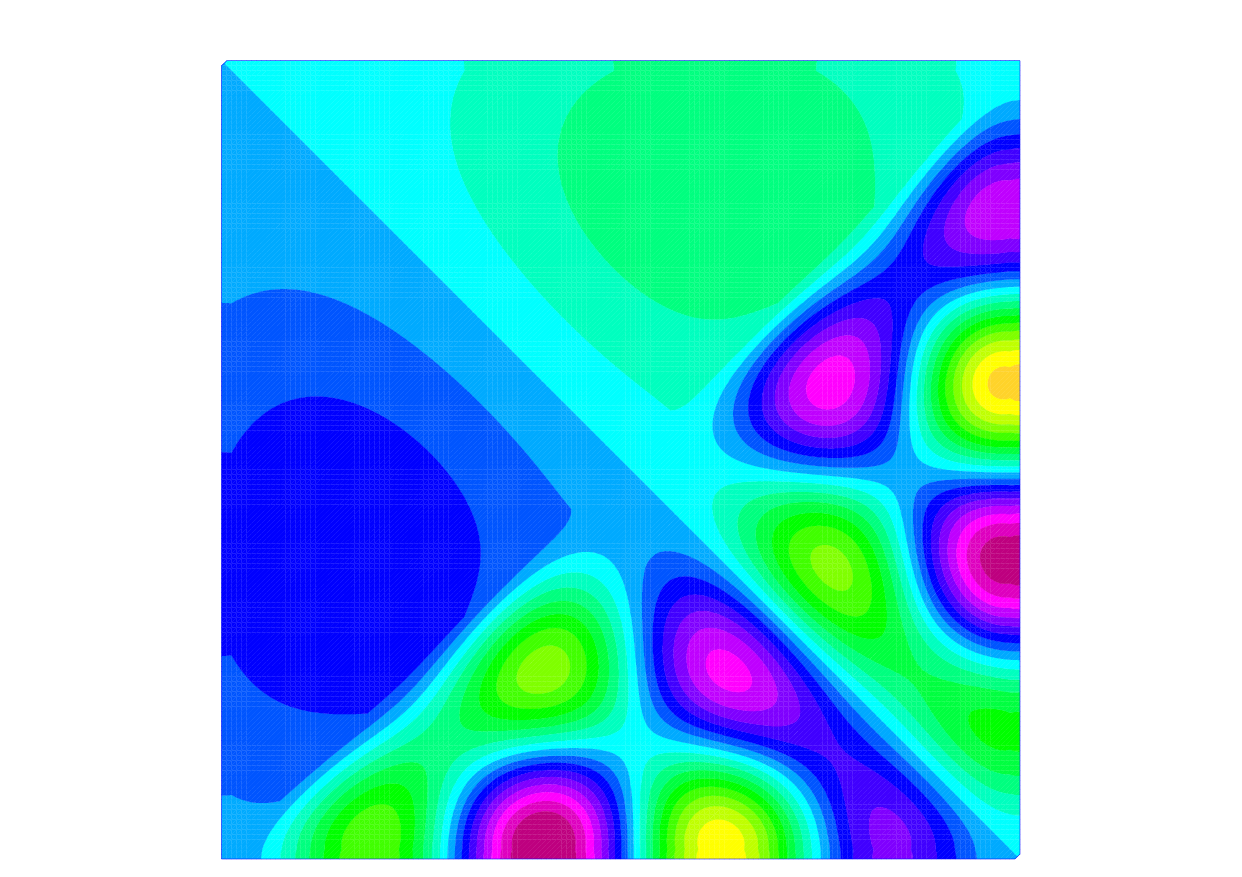} \\[-1ex]
		& $\lambda = 0.057$ & $\lambda = -0.014$ & $\lambda = -3.460$ & $\lambda = -0.052$ \\
	\end{tabu}
	\caption{Example eigenfunctions on the central subdomain when $N=25$ and $n_\text{glob} = 800$: In the first three columns, we plot qualitatively similar eigenfunctions, computed (left-to-right) by $\Delta$-GenEO, $\mathcal{H}$-GenEO when $\kappa = 1000$, and $\mathcal{H}$-GenEO when $\kappa = 10000$. This illustrates how eigenfunctions of~\eqref{eq:hg} are affected by the indefiniteness in $b_j$, relative to the size of $\kappa$. In addition, as $\kappa$ increases the $\mathcal{H}$-GenEO eigenproblem enriches the coarse space with ``wave-like'' eigenfunctions that are not seen for $\Delta$-GenEO; one of the many examples when $\kappa = 10000$ is plotted in the final column. While the top row explores the homogeneous case, the bottom row demonstrates the effect of heterogeneity in $a(\vec{x})$ for the diagonal layers problem: For $N=25$, $a(\vec{x}) = a_\text{max} = 10$ in the upper-left triangle ($x_2 > x_1$) and $a(\vec{x}) = a_\text{min} = 1$ in the lower-right triangle ($x_2 < x_1$) of the central subdomain. Note that variation in the eigenfunctions is mainly confined to the low coefficient region.}
	\label{Fig:Eigenplots}
\end{figure}

These observations align with the fact that eigenfunctions appear qualitatively similar for $\Delta$-GenEO and $\mathcal{H}$-GenEO when $\kappa$ is small. As seen in Fig.~\ref{Fig:Eigenplots}, once $\kappa$ increases the $\mathcal{H}$-GenEO eigenfunctions change: the type of eigenfunctions produced by $\Delta$-GenEO remain, albeit perturbed, but now we have further eigenfunctions which include more oscillatory behaviour in the interior of the subdomain; such features are not found with $\Delta$-GenEO where higher oscillations only appear near the boundary.

\section{Conclusions}
\label{sec:conclusions}

In this work we have summarised how the forthcoming analysis in \cite{Bootland:2021:OSM} can be applied to a GenEO-type coarse space for heterogeneous indefinite elliptic problems. We provide numerical evidence supporting these results and a comparison with a more effective GenEO-type method for highly indefinite problems but for which no theory is presently available. For mildly indefinite problems these two approaches perform similarly, providing the first theoretical insight towards explaining the good behaviour of the $\mathcal{H}$-GenEO method for challenging heterogeneous Helmholtz problems.

\bibliographystyle{spmpsci}
\bibliography{references}

\end{document}